\documentclass[12pt]{amsart}
\usepackage{amssymb,amscd,amsthm, verbatim,amsmath,color,fancyhdr, mathrsfs, bm, mathtools,soul,verbatim,listings}
\usepackage{graphicx}
\usepackage{turnstile}

\usepackage{tikz,tikz-cd}
\usepackage{array}
\usepackage{multirow}
\usepackage[normalem]{ulem}
\usepackage{amsaddr}
\usepackage{epsfig,psfrag}
\usepackage{amsfonts, amssymb}
\usepackage{enumerate}
\usepackage[all]{xy}
\usepackage[linesnumbered,ruled]{algorithm2e}
\usepackage{setspace}
\usepackage{float}
\usepackage{amsfonts, amssymb}
\usepackage{amssymb,epsfig,psfrag}
\usepackage{here}
\usepackage{pict2e}
\usepackage{hyperref}
\usepackage{setspace}
\usepackage{url}
\usepackage{here}
\usepackage{pict2e}
\usepackage{hyperref}
\usepackage{blindtext}
\usepackage{mathpazo}

\def\multiset#1#2{\ensuremath{\left(\kern-.3em\left(\genfrac{}{}{0pt}{}{#1}{#2}\right)\kern-.3em\right)}}

\hypersetup{
    colorlinks=true,
    linkcolor=blue,
    filecolor=magenta,      
    urlcolor=cyan,
}
\usepackage{babel}
\usepackage[T1]{fontenc}
\usepackage[utf8]{inputenc}
\usepackage{mwe} 
\usepackage{newfloat}
\DeclareFloatingEnvironment{argument}
\addto\captionsczech{%
}

\usepackage[letterpaper, left=2.5cm, right=2.5cm, top=2.5cm,
bottom=2.5cm,dvips]{geometry}
\newtheorem{theorem}{Theorem}[section]
\newtheorem{lemma}[theorem]{Lemma}
\newtheorem{corollary}[theorem]{Corollary}
\newtheorem{proposition}[theorem]{Proposition}
\numberwithin{equation}{section}

\theoremstyle{definition}
\newtheorem{definition}{Definition}

\newtheorem{remark}[definition]{Remark}
\newtheorem{question}[definition]{Question}

\DeclareMathOperator{\Gal}{Gal}

\DeclareMathOperator{\Aut}{Aut}

\DeclareMathOperator{\disc}{disc}
\DeclareMathOperator{\res}{res}
\DeclareMathOperator{\Res}{Res}

\DeclareMathOperator{\sgn}{sgn}

\newcommand{\Proj}{\mathbb{P}}

\newcommand{\w}{\omega}
\newcommand{\z}{\zeta}

\newcommand{\la}{\langle}
\newcommand{\ra}{\rangle}

\newcommand{\inverselimitb}

\title{Overgroups of the arboreal representation of PCF polynomials}
\author{Wayne Peng}
\address{Department of Mathematics,
        National Center University,
        Taoyuan city, Taiwan}
\email{junwen.wayne.peng@math.ncu.edu.tw}

\date{September 2021}

\begin{document}

\maketitle


\begin{abstract}
Consider a number field $K$ and a rational function $f$ of degree greater than 1 over $K$. By taking preimages of $\alpha\in K$ under successive iterates of $f$, an infinite $d$-ary tree $T_\infty$ rooted at $\alpha$ can be constructed. An edge is assigned between two preimages $x$ and $y$ if $f(x)=y$. The absolute Galois group of $K$, acting on $T_\infty$ through tree automorphisms, generates a subgroup $\Gal_f^\infty(\alpha)$ in the group of all automorphisms of $T_\infty$, $\Aut(T_\infty)$.

We have discovered a new class of natural overgroups in which the image of the Galois representation attached to a PCF polynomial must reside. Moreover, we have found that the image of the Galois representation of a new PCF polynomial is isomorphic to one of these overgroups. We also investigate the structure of these overgroups for specific maps, such as normalized dynamical Belyi polynomials, and show that the normal subgroups of these overgroups form a unique chief series. This allows us to bound the number of generators through group-theoretic analysis.
\end{abstract}


\section{Introduction}
Let $K$ be a number field and let $\alpha\in K$. Consider a polynomial $f\in K[z]$ of degree $d\geq 2$. The Galois group of $f^n(z)-\alpha$ over $K$, where $f^n$ is the $n$-th iterate of $f$, is known as the \textit{$n$-th dynamical Galois group} or the \textit{arboreal Galois group}, denoted by $\Gal_f^n(\alpha)=\Gal(f^n(z)-\alpha/K)$.

Under certain conditions, such as when neither $\alpha$ is a periodic point of $f$ or a preimage of $\alpha$ is not a branch point of $f$, the Galois group $\Gal_f^n(\alpha)$ can be embedded into the automorphism group of a $d$-ary tree $T_n(d)$ of level $n$. The vertices of this tree are made up of the inverse images of $\alpha$ under $f^{n}$ for $n=1,2,\ldots$, and an edge is drawn between two vertices $x$ and $y$ if $f(x)=y$. The point $\alpha$ is called the \textit{base point} or the \textit{root} of $T_n(d)$.

This embedding of $\Gal_f^n(\alpha)$ into the automorphism group of the tree is known as the \textit{arboreal representation}. Figure~\ref{3-ary tree} shows a tree of a cubic polynomial as an example.

We study the structure of the overgroups for PCF polynomials, including normalized dynamical Belyi polynomials. Our findings show that the normal subgroups of the overgroups form a unique chief series. This leads to the ability to bound the number of generators in terms of group-theoretic analysis. The results of our study give a new understanding of the relationship between the arboreal Galois groups and the automorphism groups of $d$-ary trees for PCF polynomials.

\begin{theorem}\label{intro-theorem-1}
Let $f$ be a degree $d$ PCF polynomial over a number field $K$, and let $\alpha\in K$ is not periodic under $f$. Let $a_f$ be the leading coefficient of $f$, let $\mathcal{C}_f$ be the set of all critical points of $f$, let $L$ be the minimal integer such that $f^L(\mathcal{C}_f)$ is a periodic set, and let $O$ be the minimal positive integer such that $f^{L+O}(\mathcal{C}_f)=f^{L}(\mathcal{C}_f)$.
\begin{enumerate}
    \item\label{subclass:1} If the degree of $f$ is odd and $L\leq 1$, then $\Gal_f^n(\alpha)$ is isomorphic to a subgroup of $E_n^{2O}(d)$.
    \item If the degree of $f$ is odd and $L > 1$, then $\Gal_f^n(\alpha)$ is isomorphic to a subgroup of $F_n^{(L+2O-1,L-1)}(d)$.
    \item If the degree of $f$ is even, then $\Gal_f^n(\alpha)$ is isomorphic to a subgroup of $E_n^{(m,m')}(d)$ with $(m,m')$ as follows:
    \begin{enumerate}
        \item If $L=0$, then $(m,m') = (O+1,1)$.
        \item If $L>1$, then $(m,m') = (L+O,L)$..
    \end{enumerate}
\end{enumerate}
\end{theorem}
It is important to note that the odd degree normalized Belyi polynomials, which have three fixed critical points $0$, $1$, and $\infty$, are a subclass of \ref{subclass:1} in Theorem~\ref{intro-theorem-1} with $L=0$ and $O=1$. The arboreal representations of normalized dynamical Belyi maps have been studied in depth in \cite{B-F-H-J-Y-2017} and \cite{BEK2020}.

The notation $E_n^m(d)$, $E_n^{(m,m')}(d)$, and $F_n^{(m,m')}(d)$ are defined in Section~\ref{sec:2.1}. They are the kernel of some sign function. A $d$-ary tree of level $n$ has $d^n$ leaves, or vertices of degree $1$. The group of automorphisms of this tree, $\Aut(T_n(d))$, can be embedded into the symmetric group $S_{d^n}$. In \cite{B-F-H-J-Y-2017} and \cite{BEK2020}, they use this embedding to define a natural sign function, denoted by $\sgn$, on $\Aut(T_n(d))$, and then define
\[
\sgn_m\coloneqq\sgn\circ\res_m
\]
where $\res_m$ is the restriction from $\Aut(T_n(d))$ to $\Aut(T_m(d))$ for any $n\geq m$. Thus, for odd $d$, we can define
\[
E_n^m(d)\coloneqq \Aut(T_n(d))\cap\ker(\sgn_m).
\]
We define a new sign function in Section~\ref{sec:2.1}, $\sgn_2^{(m,m')}$ for defining $F_n^{(m,m')}$ and $E_n^{(m,m')}$. 

The study of arboreal Galois groups has gained popularity due to its applications in number theory. Odoni \cite{Odoni1985} first demonstrated how the description of arboreal Galois groups can be used to study the density of prime divisors in certain sequences defined dynamically (see subsequent works \cite{Rafe2008,HJM2015,Looper2019} and \cite{Juul2019}). Typically, the $n$-th dynamical Galois group of a degree $d$ polynomial is isomorphic to the group of automorphisms of a $d$-ary tree $\Aut(T_n(d))$, which is equivalent to the $n$-fold wreath product $[S_d]^n$ of the symmetric group $S_d$ on $d$ elements (see \cite{juulthesis,2Odoni1985}). The fundamental question is to determine when Odoni's index $[\Aut(T_\infty(d)):\Gal_f^{\infty}(\alpha)]$ is finite. Jones (see~\cite{Rafe2013}) has proposed that the question of Odoni's index can be seen as an analogue of Serre's open image theorem (see~\cite{serre}), and many authors have provided positive answers (see~\cite{Rafe2008}, \cite{Hindes2013}, \cite{GNT2013}, \cite{BridyTucker2019}, and \cite{Li2020}).

 Jones and others (see \cite{Rafe2013,pink2013profinite,GottesmanTang2010,BHL2017,RafeMichelle2014} and \cite{benedetto2019}) have shown that the Odoni index of any PCF map is infinite, which can be seen as a dynamical analog of complex multiplication. This result is established via the exact sequence,
\[
1\to\Gal(f^n(x)-t/\overline{K})\to\Gal(f^n(x)-t/K)\to \Gal(K'/K)\to 1,
\]
where $t$ is transcendental over $K$, $\overline{K}$ is an algebraic closed field of $K$, and $K'$ is the intersection of the field extension of $f^n(x)-t$ and $\overline{K}$. The \textit{geometric Galois group} $\Gal(f^n(x)-t/\overline{K})$ is isomorphic to the \textit{iterated monodromy group} of $f$. The iterated monodromy group is a quotient group of the fundamental group $\pi_1(\Proj^1(\overline{K})\setminus f^{-1}(P_f))$. Since $P_f$ is finite, the fundamental group is finitely generated, hence the geometric Galois group is finitely generated. On the other side of the exact sequence, $\Gal(K'/K)$ is either finite or finitely generated. These two facts imply that the \textit{arithmetic Galois group} $\Gal(f^n(x)-t/K)$ is finitely generated. Therefore, any specialized $t$, i.e. $t\in K$, the Galois group is finitely generated. On the other hand, the $\infty$-fold wreath product $[S_d]^\infty$ is not finitely generated, so the arithmetic Galois group has to be isomorphic to an infinite index subgroup of $\Aut(T_\infty(d))$. Our argument and statement of Theorem~\ref{intro-theorem-1} provides a viewpoint from group theory on this question.

The normal subgroup structure of $E_n^2(d)$ is not fully understood. The next result investigates the structure of the normal subgroup of $E_n^2(d)$, which represents the dynamical Galois group of odd degree normalized dynamical Belyi maps at a carefully chosen basepoint. The normal subgroups provide a group theoretic perspective on counting the rank of $E_\infty^2(d)$, for which only an indirect proof was known. We denote the rank of a group $G$ as $d(G)$, which is the minimum cardinality of a generating set of $G$. The second result is stated in the following theorem:

\begin{theorem}\label{intro-theorem-2}
For odd $d$, the rank of the group $E_n^2(d)$ is 2 for all $n$. Additionally, $E_n^2(d)$ has a unique chief series.
\end{theorem}

The task of finding an upper bound for $d(G)$ through the degrees of subgroups or quotient groups of $G$ is difficult. In~\cite{DetomiLucchini2003}, the authors have shown that $d(G)=d(G/N)$ if and only if there are not many factors in $G/N$ that are equivalent to $N$. An alternative approach was given in~\cite{DallaVoltaLucchini}, where it was shown that $d(G)=\max\{2,d(G/N)\}$ if $G$ is not cyclic and $N$ is the unique minimal subgroup of $G$. Hence, by finding a unique minimal normal subgroup $N$ in $E_n^2$, and repeating the process for the quotient group $E_n^2/N$, it is possible to bound the rank of $E_n^2$. This process is similar to the concept of the \textit{chief series} of a finite group $G$, which is a sequence of normal subgroups $N_i$ such that $1=N_1\triangleleft N_2\triangleleft\cdots\triangleleft N_k=G$ and each chief factor $N_{i+1}/N_{i}$ is a minimal normal subgroup of the quotient group $G/N_i$ (see~\cite{Isaacs2009} for more details). By applying this analysis, the rank of large groups can be found and the structure of their normal subgroups can be clearly defined.

The determination of the dynamical Galois group of a polynomial is often a challenging task. However, we can use Theorem~\ref{intro-theorem-1} as a criterion to search for new PCF polynomials. It should be noted that the critical portrait, defined as $0\mapsto 1$, $1\mapsto 0$, and $\infty\mapsto\infty$, had not been studied previously. Thus, we present it as our final theorem.
\begin{theorem}\label{intro-theorem-3}
Let $K$ be a number field, and let $f(z)=2z^3-3z^2+1$. If there exist primes $\mathrm{p}$ and $\mathrm{q}$ of $K$ lying above $2$ and $3$, respectively, and an element $\alpha\in K$ such that $v_\mathrm{q}(\alpha)=1$ or $v_\mathrm{q}(1-\alpha)=1$, and $v_\mathrm{p}(\alpha)=1$ or $v_\mathrm{p}(1-\alpha)=1$, then the following holds for every $n\geq 1$:
\begin{enumerate}
    \item The polynomial $f^n$ is irreducible over $K$.
    \item The $n$-th dynamical Galois group, $\Gal_f^n(\alpha)$, is isomorphic to $E_n^2(3)$.
\end{enumerate}
\end{theorem}

The structure of this paper is as follows: In Section 2, we introduce our overgroups, provide a brief overview of the relationship between the abelianization and the wreath product, and calculate the discriminant of $f^n-\alpha$. Most of the content in Section 2 is well-known. In Section 3, we present our main results. Finally, in Section 4, we address open questions and discuss the challenges of generalizing the arguments in this paper for more general cases.

Throughout this paper, we use the following symbols:
\begin{itemize}
    \item $1$ is used to represent either the integer $1$, an identity map, or the identity of a group. If $1$ is used to denote an identity, the group will be specified in the context.
    \item $d$, $n$, and $m$, with or without sub-indices, represent integers.
    \item The symbol $\ast$ is used to denote any arbitrary element. The set to which the element belongs will be specified in the context.
    \item For $\sigma,\tau\in G$, we denote $\sigma^\tau$ for the conjugacy $\tau\sigma\tau^{-1}$.
    \item $A_d$ represents an alternating group of degree $d$, and $S_d$ represents a symmetric group of degree $d$.
    \item The set of all finite critical points of $f$ is denoted by $\mathcal{C}_f\coloneqq (f')^{-1}(0)$, and $f(\mathcal{C}_f)=\{f(c)\mid c\in\mathcal{C}_f\}$.
    \item The splitting field of $f^n-\alpha$ over $K$ is denoted by $K_f^n(\alpha)$, and its corresponding Galois group is denoted by $\Gal_f^n(\alpha)$.
    \item An object $\alpha$, which can be a point or a set, is said to be preperiodic under a map $f$ if $\{f^n(\alpha)\mid n=1,2,\ldots\}$ is a finite set. If $\alpha$ is periodic, then $f^n(\alpha)=\alpha$ for some integer $n$. The tail length of a preperiodic object is the minimum integer such that $f^n(\alpha)$ becomes a periodic object. A preperiodic object is said to be strictly preperiodic if its tail length is not zero.
    \item $f^{-n}(\alpha)$ represents the set of all $n$-th preimages of $\alpha$ under $f$. An $n$-th preimage of $\alpha$ under $f$ is a root of $f^n(z)-\alpha$, and multiple roots are considered as different elements in the set.
\end{itemize}


\section{Preliminaries}\label{sec:2}
\subsection{Wreath product and arboreal representation}\label{sec:2.1}

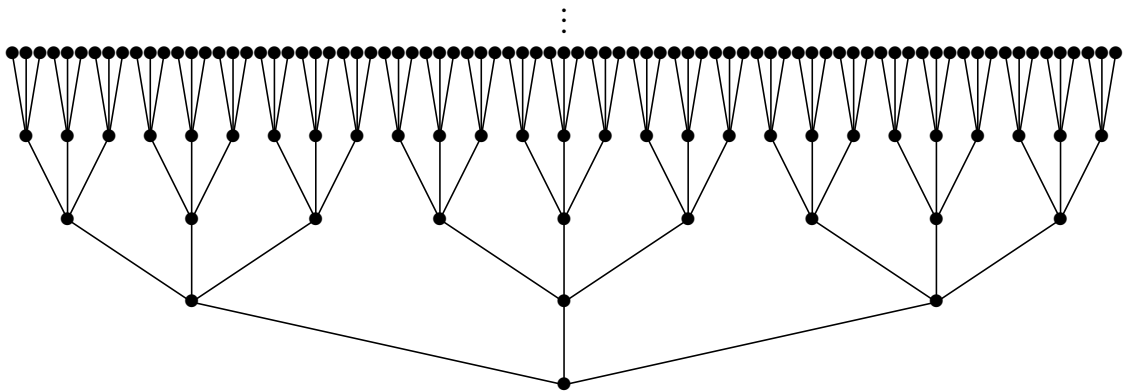
\begin{figure}[ht]
\centering
\begin{tikzpicture}
[scale = 0.55]
	\node (top) at (8.5,9) {$\vdots$};
	\node (0) at (8.5,0) {$\color{black}\bullet$};
	\foreach \x in {40,41,...,120}
		\node (\x) at (1/3*\x -80/3+8.5,8) [inner sep=0.5pt] {$\color{black}\bullet$};
	\foreach \x in {13,...,39}
		\node (\x) at (\x - 26 +8.5,6) [inner sep=3pt] {$\color{black}\bullet$};
	\foreach \x in {4,5,...,12}
		\node (\x) at (3*\x-24 +8.5,4) [inner sep=3pt] {$\color{black}\bullet$};
	\foreach \x in {1,2,3}
		\node(\x) at (9*\x - 18+8.5,2) [inner sep=3pt]  {$\color{black}\bullet$};
	\foreach \x [evaluate=\x as \y using int(3*\x + 1),evaluate=\x as \z using int(3*\x + 2),evaluate=\x as \w using int(3*\x + 3)] in {0,1,...,39}{
		\draw[line width=.6pt] (\x.center) -- (\y.center);
		\draw[line width=.6pt] (\x.center) -- (\z.center);
		\draw[line width=.6pt] (\x.center) -- (\w.center);
	}
\end{tikzpicture}
\caption{A 3-ary tree of 4 levels}\label{3-ary tree}
\end{figure}

Let $T_n(d)$ be a $d$-ary tree with $n$ levels (refer to Figure~\ref{3-ary tree} for an example). It has $d^n$ leaves and a total of $1 + d + \cdots + d^n$ vertices. The level of a vertex on the tree is defined as the distance from the vertex to the root. For $m \leq n$, the $m$-th level of $T_n(d)$ is the set of vertices with level $m$.

Our arguments in this paper rely on a specific labeling of the tree $T_n(d)$. To ensure rigor, we will define this labeling explicitly, but will not elaborate further on it in the rest of the paper.

There are two types of labelings for $T_n(d)$. The \textit{regular labeling} assigns a tuple $(l_m,\ldots, l_1)$ in $\{0,1,2,\ldots,d^m-1\}$ to each vertex at level $m$. The vertex with label $(l_m,\ldots, l_1)$ is at a distance of 1 from the vertex with label $(l_{m-1},\ldots, l_1)$. The root of the tree is labeled with an empty tuple $()$. This labeling procedure provides a recursive relationship for labeling all vertices in the tree.

We can identify some specific subtrees using the regular labeling. For example, for each tuple $(l_{m},\ldots,l_{1})$ with $m\leq n$, the set of vertices $(\ast,\cdots,\ast,l_{m},\ldots,l_1)$ in $T_n(d)$ form a subtree that is isomorphic to $T_{n-m}(d)$, rooted at $(l_m,\ldots, l_1)$. This subtree is referred to as the \textit{level $n-m$ subtree} and is denoted by $T_{(l_m,\ldots, l_1)}(d)$.

To simplify the labeling system, we introduce a new labeling system specifically for the vertices on the $n$-th level of $T_n(d)$. We denote this labeling system as $I(T_n)={0,1,2,\cdots, d^{n}-1}$, where each vertex $(l_n,\ldots, l_1)$ is labeled by the following expression:
\[
(l_n-1)d^{n-1}+(l_{n-1}-1)d^{n-2}+\cdots+(l_{1}-1),
\]
Two indices $i$ and $j$ belong to the same level $m$ subtree $T_{(l_m,\ldots, l_1)}(d)$ if $i\equiv j\ (\text{mod}\ d^{m-1})$. To simplify further, we group indices into blocks based on the subtree they belong to. The \textit{block} $\mathcal{B}(T_n/T_m)(l_{m},\ldots, l_1)$ consists of indices that belong to the same level $m$ subtree rooted at $(l_{m},\ldots, l_1)$. We use $\mathcal{B}(T_n/T_m)$ to denote the set of all blocks. An index $i$ belongs to a level $m$ subtree $T_{(l_m,\ldots, l_1)}$ if it is in the corresponding block $\mathcal{B}(T_n/T_m)(l_m,\ldots, l_1)$.

We will not go into the explicit definition of $\Aut(T_n(d))$ as it can be found in \cite{B-F-H-J-Y-2017}. For readers unfamiliar with this concept, we suggest looking it up in the reference for a more detailed explanation. The labelings of $T_n(d)$ induce an embedding $[S_d]^n\hookrightarrow S_{\# I(T_n(d))}$, which connects the symmetric group theory to the labelings of the tree. For example, an index $i\in I(T_n(d))$ is considered to be in the "support" of $\sigma\in\Aut(T_n)$ if $\sigma(i)\neq i$.

We will adopt the definition of $\sgn$ from \cite{B-F-H-J-Y-2017} and use it to construct the groups $E_n^m$, $E_n^{(m,m')}(d)$, and $F_n^{(m,m')}(d)$. The details will be given after we briefly discuss the wreath product.

The wreath product of $H$ by $G$ is a group extension of $G$ by $H^{|I|}$, where $H$ is an arbitrary group and $G$ is a group that acts on a set $I={1,2,\ldots, d-1}$. An element in the wreath product $H\wr_{I} G$ is expressed as
\[
((h_1,\ldots, h_d);g)
\]
 with $h_i\in H$ and $g\in G$. The operation is defined as 
\[
((h_i)_{i\in I}; g)((h'_i)_{i\in I} ; g')=((h_{g'(i)})_{i\in I}(h'_i)_{i\in I};gg')
\]
where $(h_{g'(i)})_{i\in I}(h'_{i})_{i\in I}$ is simply the component-wise product. The following simple decomposition of an element in a wreath product is frequently used in the argument:
\[
((h_1,\ldots, h_d);g)=((1,\ldots, 1);g)
((h_1,\ldots, h_d);1)
\]
For clarity, the index set $I$ is omitted when it is clear from the context. It is important to note that for $n>m$, $\Aut(T_n(d))\cong \Aut(T_m(d))\wr_{I(T_{n-m}(d))}\Aut(T_{n-m}(d))$. Elements in $\Aut(T_m(d))\wr\Aut(T_{n-m}(d))$ are expressed as
\[
((a_i)_{i\in I(T_{n-m}(d))};b)
\]
when they are considered as part of the group $G\subseteq \Aut(T_m(d))\wr\Aut(T_{n-m}(d))$. The projection $\res:H\wr G\to G$ is naturally a homomorphism, defined as $\res((h_i);g)=g$. Combined with the above observation, for any $n\geq m$, the restriction
\[
\res_m:\Aut(T_{n-m}(d))\wr\Aut(T_m(d)) \to \Aut(T_m(d)),
\]
is defined. Note that $\res_n$ is the identity map. The restriction $\res_m$ is referred to as the \textit{restriction to the level $m$} or simply as the \textit{restriction}.

Let's recall the sign of a permutation. We define the sign of a permutation $\sigma$ as $+1$ if $\sigma$ is even and $-1$ if $\sigma$ is odd. To define a sign on the automorphisms of a $d$-level $n$-ary tree, we label the leaves of the tree, enabling us to embed $\Aut(T_n(d))$ into $S_{d^n}$. The composition of the embedding with the natural sign of permutations defines a sign on $\Aut(T_n(d))$. Two different labelings result in different embeddings that are conjugated by an element in $S_{d^n}$, and conjugation doesn't change the sign, so we have a well-defined homomorphism $\sgn:\Aut(T_n(d))\to C_2={\pm 1}$. Finally, for a $d$-level $n$-ary tree with $n\geq m$, we define the \textit{natural sign}
\[
\sgn_{m}=\sgn\circ \res_{m}
\]
and refer to $\sgn_m$ as the \textit{sign restricted to level $m$} or simply the \textit{restricted sign}. Here are some propositions about the sign on the automorphisms of a $d$-level $n$-ary tree.
\begin{lemma}\label{lemma:sgn}
Let $T_n(d)$ be a $d$-ary tree of $n$ levels, with $m>n$, and let $\sigma=((b_i){i\in I(T_m(d))};a)\in \Aut(T_n(d))\cong \Aut(T_{n-m}(d))\wr\Aut(T_m(d))$. Then, we have
\[
\sgn(\sigma)=
\begin{cases}
    \sgn(a)\prod_{i\in I(T_m(d))}\sgn(b_i), &\text{ for odd }d,\\
    \prod_{i\in I(T_{m}(d))}\sgn(b_i), &\text{ for even }d.
\end{cases}
\]
\end{lemma}
\begin{proof}
First, we note that 
\[
\sgn(((b_i)_{i\in I};a))=\sgn(((1)_{i\in I}, a)))\sgn(((b_i)_{i\in I};1)))\in\Aut(T_{n-m}(d))\wr\Aut(T_n(d)).
\]
Since $b_i$ acts disjointly on the tree $T_n(d)$, it is clear that $\sgn(((b_i)_{i\in I};1))=\prod_{i\in I}\sgn(b_i)$. Now, we focus on $\sgn(((1)_{i\in I}; a))$. Let us express $a$  as a product of disjoint cycles $a_1\cdots a_k$ for some integer $k$, i.e. $((1);a)=((1);a_1)\cdots((1);a_k)$. Each permutation $((1);a_i)$ creates $d^{n-m}$ many disjoint cycles in the same length as $a_i$ on the $n$-th level of the tree $T_n(d)$, and the cycles associated with $a_i$ are disjoint from the cycles associated with $a_j$ for distinct $i$ and $j$. Therefore, the sign of $((1)_{i\in I}, a)$ is
\[
\prod_{i=1}^k\sgn(a_i)^{d^{m-n}}=
\begin{cases}
1 \text{, if }d\text{ is even;}\\
\prod_{i=1}^k\sgn(a_i)\text{, otherwise.}
\end{cases}
\]
\end{proof}

For $0<m<n$, we can define an alternative sign function
\begin{align*}\sgn^m:\Aut(T_n(d))\cong\Aut(T_{n-m}(d))\wr\Aut(T_m(d))&\to \{\pm 1\}\\
((b_i)_{i\in I}; a)\mapsto\prod_{i\in I}\sgn(b_i)
\end{align*}
with $b_i\in \Aut(T_{n-m}(d))$ and $a\in\Aut(T_{m}(d))$. It is an easy exercise to check $\sgn^m$ is a group homomorphism.

For $n\geq m>m'>0$, we define the following to construct an overgroup of $\Gal_f^n(\alpha)$ for a PCF map $f$. Let
\begin{equation}\label{eq:sgn3}
\sgn^{(m,m')}_1=(\sgn^{m'}\circ\res_{m})\cdot\sgn_{m'}
\end{equation}
and
\begin{equation}\label{eq:sgn2}
\sgn^{(m,m')}_2=\sgn^{m'}\circ \res_{m}.
\end{equation}
\begin{remark}
The integer $n$ is included to define the sign functions on $Aut(T_n(d))$, but the functions themselves only depend on the tree structure up to level $m$. Thus, we suppress $n$ in the notation. It is important to note that for even $d$, the natural sign that we defined above only takes into account the level $1$ subtrees on the top of the tree, and the permutation on the lower branches is not considered. As a result, $\sgn_1^{(m,m')}$ is not a natural sign for even $d$ and is a natural sign for odd $d$. The similar conclusion can be drawn for $\sgn_2^{(m,m')}$. To clarify this discussion, the following table is provided.
\begin{center}
    \begin{tabular}{c|c|c}
    $\sgn$\textbackslash d& even & odd \\\hline
     $\sgn^{(m,m')}_1$ & & natural sign \\\hline
     $\sgn^{(m,m')}_2$ & natural sign &
    \end{tabular}
\end{center}
\end{remark}
\begin{lemma}\label{lemma:sgn2}
For $\sigma\in\Aut(T_n(d))$ and $m>m'>0$, the equation
\[
\sgn^{(m,m')}_i(\sigma)=\sgn_{m}(\sigma)\cdot \sgn_{m'}(\sigma)
\]
holds for $i=1$ when $d$ is even, or for $i=2$ when $d$ is odd. Otherwise, we have
\[
\sgn_i^{(m,m')}=\sgn_m(\sigma).
\]
\end{lemma}
\begin{proof}
The second assertion follows directly from Lemma~\ref{lemma:sgn}, so we will focus on proving the first assertion. We use the notations that we defined above, and the proposition follows from direct computation. For $i=1$ and even $d$, the equation follows from the definition. For $i=2$, we let $\sigma\in\Aut(T_n(d))\cong \Aut(T_{n-m}(d))\wr\Aut(T_m(d))$ with $n\geq m$ and note that $\sgn_{m}(\sigma)=\sgn(a)$, and so it follows that
\[
\sgn(\sigma)\sgn_{m}(\sigma)=\sgn(a)\prod_{i\in I}\sgn(b_i)\sgn_{m}(\sigma)=\prod_{i\in I}\sgn(b_i)=\sgn^m(\sigma).
\]
By Equation~\eqref{eq:sgn2}, we have
\begin{align*}
    \sgn^{(m,m')}_2&=\sgn^{m'}\circ\res_{m}\\
    &=(\sgn\cdot\sgn_{m'})\circ\res_{m}\\
    &=(\sgn\circ\res_{m})\cdot(\sgn_{m'}\circ\res_{m})\\
    &=\sgn_{m}\cdot(\sgn\circ\res_{m'}\circ\res_{m})=\sgn_{m}\cdot\sgn_{m'}.
\end{align*}
\end{proof}
We are now ready to define two subgroups of $\Aut(T_n(d))$. Denote $I(T_1(d))$ by $I$. For $n>0$ and $m> m'>0$, the groups are defined as follows:
\[
E_n^{(m,m')}(d)=
\begin{cases}
\Aut(T_n(d))\text{, if }n<m;\\
E_{n-1}^{(m,m')}(d)\wr_{I}\Aut(T_1(d))\cap \ker(\sgn_1^{(m,m')})\text{, if }n\geq m
\end{cases}
\]
and
\[
F_n^{(m,m')}(d)=
\begin{cases}
\Aut(T_n(d))\text{, if }n<m;\\
F_{n-1}^{(m,m')}(d)\wr_{I}\Aut(T_1(d))\cap \ker(\sgn^{(m,m')}_2)\text{, if }n\geq m.
\end{cases}.
\]
The $(d)$ will be omitted if the context is clear.

The domain of the sign functions in the above definition changes as $n$ increases. We embed $E_{n-1}^{(m,m')}\wr \Aut(T_1)$ and $F_{n-1}^{(m,m')}\wr\Aut(T_1)$ into $\Aut(T_n)$, which makes the $\sgn_1^{(m,m')}$ and $\sgn_2^{(m,m')}$ functions well-defined. Since $\sgn_1^{(m,m')}$ and $\sgn_2^{(m,m')}$ are the natural sign function $\sgn_{m}$ for even $d$ and odd $d$ respectively, we will simplify the notations for $E_n^{(m,m')}$ to $E_n^{m}$ and $F_n^{(m,m')}$ to $F_n^{m}$ for even and odd $d$ respectively.

The groups $E_n^{(m,m')}$ and $F_n^{(m,m')}$ serve as our target overgroups. It is natural to determine their orders, which we can calculate using the following proposition. The order of these groups is independent of the parameter $m'$, as they are kernels of certain sign functions on $\Aut(T_m(d))$. However, the parameter $m'$ does influence the structural properties of these groups.
\begin{proposition}\label{prop:unboundedindex}
We have
\[
|E_n^{(m,m')}(d)|=|F_n^{(m,m')}(d)|=
\dfrac{(d!)^{\frac{d^n-1}{d-1}}}{2^{\frac{C(d,n)}{d-1}}}
\]
where
\[
C(d,n)=\begin{cases}
0\text{, if }n<m;\\
\dfrac{d^{n-m+1}-1}{d-1}, \text{ otherwise.}
\end{cases}
\]
\end{proposition}
\begin{proof}
For $n\geq m$, we claim that the homomorphisms
\[
\begin{tikzcd}
\phi_1^{[n]}: E_{n-1}^{(m,m')}\wr \Aut(T_1(d))\arrow[r,hook]& \Aut(T_n(d)) \arrow[r,"\sgn_1^{(m,m')}"]&\{\pm 1\}
\end{tikzcd}
\]
and
\[
\begin{tikzcd}
\phi_2^{[n]}: F_{n-1}^{(m,m')}\wr \Aut(T_1(d))\arrow[r,hook]& \Aut(T_n(d)) \arrow[r,"\sgn_2^{(m,m')}"]&\{\pm 1\}
\end{tikzcd}
\]
are surjective.

When $d$ is odd, showing the claim for $\phi_1^{[n]}$ is fairly straightforward. Taking $g=((1)_{i\in I},\sigma)\in E_{n-1}\wr \Aut(T_1)$ where $\sigma$ is an odd permutation, we have $\sgn_m(g)=\sgn(\sigma)=-1$.

When $d$ is even, obviously $\phi_1^{[m]}$ is surjective because the domain is a full wreath product of $S_d$. To show the claim for $n>m$, we apply the following induction on two statements:
\begin{itemize}
    \item[$\mathcal{S}_1^{[n]}$:] The homomorphism $\varphi_n:E_{n-1}^{(m,m')}\wr \Aut(T_1)\to C_2^{m}$ by defining
    \[
    \varphi_n: a\mapsto (\sgn_{1}(a),\ldots,\sgn_{m}(a))
    \]
    is surjective.
    \item[$\mathcal{S}_2^{[n]}$:] $\varphi_n(\ker\phi_1^{[n]})\subseteq C_2^{m}$ is the codimensional $1$ variety $X_{m}X_{m'}=1$. 
\end{itemize}
We will show $\mathcal{S}_1^{[n]}\Rightarrow \mathcal{S}_2^{[n]}$ and $\mathcal{S}_2^{[n]}\Rightarrow \mathcal{S}_1^{[n+1]}$. For $\mathcal{S}_1^{[n]}\Rightarrow \mathcal{S}_2^{[n]}$, we define a homomorphism
\begin{align*}
    \psi:C_2^{m}&\to C_2\\
    (c_{1},\ldots,c_{m})&\mapsto c_{m}c_{m'},
\end{align*}
and note that $\phi_1^{[n]}=\psi\circ\varphi_n$ by Lemma~\ref{lemma:sgn2}. Since $\varphi_n$ is surjective and $\ker(\phi_1^{[n]})$ has index two in $E_{n-1}^{(m,m')}\wr\Aut(T_1(d))$, its image $\varphi_n(\ker(\phi_1^{[n]}))$ must be the variety $X_mX_{m'}=1$ in $C_2^m$, as this is the unique subgroup of index two in $C_2^m$ containing the kernel of $\psi$. For $\mathcal{S}_2^{[n]}\implies \mathcal{S}_1^{[n+1]}$, we observe that, for $b=((a_i);\sigma)\in E_n^{(m,m')}\wr\Aut(T_1)$, where $a_i=1$ for all $i\neq 1$,

\begin{align*}
\varphi_{n+1}(b) &= (\sgn_1(b),\sgn_2(b),\sgn_3(b),\ldots,\sgn_m(b))\\
&=(\sgn(\sigma),\sgn_{1}(a_1),\sgn_2(a_1),\ldots,\sgn_{m-1}(a_1)),
\end{align*}
where the last equality follows from $n$ being even and Lemma~\ref{lemma:sgn}. Since $\sgn_m(a_1)$ is not one of the components, there is no restriction on $\sgn_{m'}(a_1)$. This proves that $\varphi_{n+1}$ is surjective. Finally, since $\varphi_n$ is surjective for all $n>m$, $\phi_{1}^{[n]}$ is also surjective for all $n>m$.

For odd $d$, the argument of showing that $\phi_2^{[n]}$ satisfies the claim is almost identical to the above argument. For even $d$, we only need to change $\mathcal{S}_2^{[n]}$ to the following
\begin{itemize}
    \item[$\mathcal{S}_2^{[n]}$:] $\varphi_n(\ker\phi_2^{[n]})\subseteq C_2^{m}$ is the codimensional $1$ variety $X_{m}=1$,
\end{itemize}
and the rest of the argument is similar.

Now, we compute the orders of these groups. It is equal to $(d!)^{\frac{d^n-1}{d-1}}$ for all $n<m$ since $\Aut(T_n)\cong [S_d]^n$. We proceed by induction. Suppose the formula is true for all $n\geq m-1$, then
\[
|E_n^m|=\dfrac{1}{2}|E_n^m\wr\Aut(T_1)|=\dfrac{1}{2}\left(\dfrac{(d!)^{\frac{d^n-1}{d-1}}}{2^{\frac{d^{n-m+1}-1}{d-1}}}\right)^dd!=\left( \dfrac{(d!)^{\frac{d^{n+1}-1}{d-1}}}{2^{\frac{d^{n-m+2}-1}{d-1}}}\right).
\]
A similar argument holds for $F_n^{(m,m')}$.
\end{proof}
\begin{lemma}\label{lemma:transitive}
For odd values of $d > 2$, $E_n^m(d)$ acts transitively on the index set $I(T_n)$. This means that for any two distinct indices, there exists an element $\sigma \in E_n^m(d)$ that transposes $i$ and $j$.
\end{lemma}
\begin{proof}
We prove the statement by induction. We observe that $E_n^1$ is a subset of $E_n^m$ for any $m$, so it suffices to show that the statement holds for $E_n^1$ for any $n>0$. $E_n^1$ is isomorphic to the $n$-fold wreath product of $A_d$, which makes the statement straightforward to verify.

The base case is trivial. We assume the statement is true for some $n$. To prove it for $n+1$, if $i$ and $j$ are on the same level $k$ subtree with $k<n$, the statement is true by induction and the fact that the group is a self-similar group. If $i$ and $j$ are not in the same subtree of level $k<n$, we can find a transposition $\tau\in E_{n-k}^{1}$ that transposes the subtree of $i$ and $j$, a transposition $\sigma_i\in E_k^1$ on the subtree of $i$, that takes $\tau(j)$ to $i$, and a transposition $\sigma_j\in E_K^1$ on the subtree of $j$ that takes $\tau(i)$ to $j$. Then, we consider 
\[
    \sigma=((1,\ldots, 1,\sigma_i,1\ldots,1,\sigma_j,1,\ldots,1),\tau)
\]
where $\sigma_i$ and $\sigma_j$ are at the coordinates corresponding the subtree of $i$ and $j$ respectively. It can be verified that $\sigma$ is in $E_n^1$ and transposes $i$ and $j$."
\end{proof}
\subsection{Abelianization}
In this section, we present a well-known result regarding the relationship between the abelianization and the wreath product of a group. We provide the statement and proof for the sake of completeness. The \emph{abelianization} of a group $G$, denoted by $G^{ab}$, is defined as the largest abelian quotient group of $G$. Let $[H,K]$ be the subgroup generated by $hkh^{-1}k^{-1}$ for all $h\in H$ and $k\in K$. The commutator subgroup of $G$ is $[G,G]$, and $G^{ab}\cong G/[G,G]$.

The following theorem deduces a relation between the wreath product and the abelianization. 
\begin{lemma}
Let $G$ be the semidirect product of $H$ by $N$. Then, $G^{ab}=(N\rtimes H)^{ab}=(N^{ab})_H\times (H^{ab})$ where the subscription $H$ denotes taking the coinvariants with respect to $H$.
\end{lemma}
\begin{proof}
We treat $N$ and $H$ as subgroups of $G$. Since $[G,G]=\la [N,N],[N,H],[H,H]\ra$, and $N\cap H$ is trivial, we have 
\begin{align*}
    G^{ab}&=\dfrac{G}{[G,G]}=\dfrac{N\rtimes H}{\la [N,N],[N,H],[H,H]\ra}\\
    &\cong \dfrac{(N^{ab})\rtimes H}{\la [H,N^{ab}],[H,H]\ra}\cong \dfrac{(N^{ab})_H\times H}{\la [H,H]\ra}\\
    &\cong (N^{ab})_H\times H^{ab}.
\end{align*}
For the second isomorphism from the back, we note that $[H,N^{ab}]\subseteq N^{ab}$, where $H$ is actually the quotient $H/H\cap [N,N]$, so $[H,N^{ab}]\cap H=\{1\}$. Moreover $n$ and $h^{-1}nh$ are on the same orbit of $n$ by the group action $N$ on $H$, so $N^{ab}/[H,N^{ab}]$ is isomorphic to $(N^{ab})_{H}$. 
\end{proof}
Using this lemma, we can compute the abelianization of a wreath product.
\begin{lemma}\label{lemma:abelinaiizationofwreathproduct}
Let $G$ and $H$ be two groups and let $G$ act faithfully and transitively on a finite index set $I$. Then, the abelianization of the wreath product $H\wr_I G$ is equal to the direct product of the abelianizations of $G$ and $H$, that is, $(H\wr_I G)^{ab} = G^{ab} \times H^{ab}$.
\end{lemma}
\begin{proof}
Let $B=H^r$, where $G$ acts faithfully and transitively on the index set $I={1,2,\ldots,r}$ of $B$. The lemma states that $G\wr H=G\ltimes B$, which implies that $(G\wr H)^{ab}=G^{ab}\times (B^{ab})_G$. Furthermore, $B^{ab}=(H^r)^{ab}=(H^{ab})^r$. To show that $((H^{ab})^r)_G=H^{ab}$, we define the following map and prove that it is a surjective homomorphism with $[G,B]$ as its kernel.

Let $\phi: (H^{ab})^r\to H^{ab}$ be defined by
\[
\phi((h_i)_{i})=\prod_{i=1}^r h_i.
\]
Since $H^{ab}$ is abelian, the map is a homomorphism. It is also surjective, as we can choose $(h_i)$ such that $h_1$ is any element in $H^{ab}$ and $h_i=1$ for all $i\neq 1$. The kernel of $\phi$ contains $[(H^{ab})^r,G]$. Conversely, if $\phi((b_i))=\sum_{i=1}^rb_i=0$, then
\[
b_1^{-1}=\prod_{i=2}^r b_i.
\]
 there exists a $g_i$ such that $g_i(i)=1$ for each $i$. For $i\neq 1$, let us define $\mathbf{b}{i_0}=(b{i}')$, where $b_{i}'=1$ for all $i\neq i_0$ and $b'{i{0}}=b_{i_{0}}$. Thus, $g_{i_0}^{-1}\mathbf{b}^{-1}{i_0}g{i_0}=(b''{i})$ with $b_1''=b'{i_0}$ and $b_i=1$ for all $i=2,3,\ldots,r$, and$g_{i_0}^{-1}\mathbf{b}^{-1}_{i_0}g_{i_0}=(b''_{i})$ with $b_1''=b'_{i_0}$ and $b_i=1$ for all $i=2,3,\ldots,r$ and
\[
(b_i)=\prod_{i_0=2}^rg^{-1}_{i_0}\mathbf{b}^{-1}_{i_0}g_{i_0}\mathbf{b}_{i_0}.
\]
Therefore, $(b_i)\in [(H^{ab})^r,G]$.
\end{proof}

\subsection{Discriminant}
In this section, we aim to demonstrate that the discriminant of the $n$-th iteration of a PCF polynomial $f$ is a square in a fixed field when $n$ is large enough. The formulas for this calculation have been given in \cite{Isaacs2009,Rafe2013} and \cite{CH2012}. However, we provide a refined formula by more precise calculation of the sign. The formula is:
\begin{proposition}\label{prop:discriminant}
Let $f$ be a polynomial over a field $K$, let $a_f$ be its leading coefficient, let $d$ be its degree, let $\alpha\in K$, and let $n\geq 2$ be an integer. The discriminant of the $n$-th iteration of $f$, with multiple roots considered to be different elements, is given by:
\begin{equation}\label{eq:disc}
\disc(f^n(x)-\alpha)=(-1)^{A(d,n)}a_f^{B(d,n)}d^{d^n}\disc(f^{n-1}(x)-\alpha)^{d}\prod_{c\in \mathcal{C}_f}(f^n(c)-\alpha),
\end{equation}
where
\begin{align*}
A(d,n) & =d^n\left(\dfrac{d^n-1}{2}+\dfrac{d^{n-1}-1}{2}\right), \\
B(d,n) & =d^{2n-1}-1.
\end{align*}
In particular, $(-1)^{A(d,n)}=\left(\frac{d}{4}\right)$ is equal to 1 if $d$ is even or if $d \equiv 1 \mod 4$, and equal to -1 otherwise.
\end{proposition}
\begin{proof}
We will focus on calculating the power of $-1$ in the discriminant. Let $c_1, c_2, \ldots$ be constants that only involve the leading coefficients of $f$ and the degree of $f$. Using the definition of the determinant, we have
\begin{align*}
\disc(f^n(x)-\alpha) & =(-1)^{d^n(d^n-1)/2}c_1\prod_{a\in f^{-n}(\alpha)}(f^n)'(a)\\
& =(-1)^{d^n(d^n-1)/2}c_1\prod_{a\in f^{-n}(\alpha)}(f')(f^{n-1}(a))f'(f^{n-2}(a))\cdots f'(f(a))f'(a).
\end{align*}
Note that
\begin{align*}
\prod_{f^{n}(a)=\alpha}&(f')(f^{n-1}(a))f'(f^{n-2}(a))\cdots f'(f(a))\\&=\left(\prod_{f^{n-1}(a)=\alpha}(f')(f^{n-2}(a))f'(f^{n-3}(a))\cdots f'(a)\right)^{d}\\
&=\left((-1)^{d^{n-1}(d^{n-1}-1)/2}c_2\disc(f^{n-1}(x)-\alpha)\right)^d
\end{align*}
The product in the tail of the equation above is $\Res(f^n(x)-\alpha,f'(x))/c_3$, where $\Res(P,Q)$ is the resultant of the polynomials $P$ and $Q$. To obtain Equation~\ref{eq:disc}, we use the formula
\[
\Res(P,Q)=(-1)^{\deg P\deg Q}\Res(Q,P).
\]
Since either $\deg(f^n(x)-\alpha)$ or $\deg(f'(x))$ is even, the product does not contribute any additional $-1$, giving us the desired result.
\end{proof}
We now present a straightforward application of Lemma~\ref{lemma:abelinaiizationofwreathproduct}.
\begin{proposition}\label{prop:quadraticextension}
If $\Gal_f^n(\alpha)$ is isomorphic to $[S_d]^n$, then $L=K({\sqrt{\disc(f^m(x)-\alpha)}\mid m=1,\ldots, n})$ has index $2^{n}$ over $K$, and any quadratic subextension of $K_f^{n}(\alpha)$ is contained in $L$.
\end{proposition}
\begin{proof}
By Lemma~\ref{lemma:abelinaiizationofwreathproduct}, the abelianization $[S_d]^{ab}$ is the $n$-fold direct product of $S_d^{ab}=C_2$. Therefore, if we can find an abelian extension $K^{ab}$ of $K$ contained in $K_f^n(\alpha)$ with the index $[K^{ab}:K]=2^n$ that is built from adjoining square roots, then we find all possible quadratic extensions of $K$ contained in $K_f^{n}(\alpha)$.

We will show that $K(\{\sqrt{\disc(f^m(x)-\alpha)}\mid m=1,\ldots,n\})$ is a subfield of $K_f^n(\alpha)$ that is of degree $2^n$ over $K$ by showing that $K(\sqrt{\disc(f^{n}(x)-\alpha)})$ is not in $K(\{\sqrt{\disc(f^m(x)-\alpha)}\mid m=1,\ldots,n-1\})$ for all $m=1,\ldots,n-1$. Note that this assertion naturally implies $\disc(f^m(x)-\alpha)$ is not in $K$.

Note that $(((12),1,\ldots,1);1)\in \Gal(f^{n}-\alpha)\cong \Aut(T_1)\wr\Aut(T_n)$, and the automorphism swaps $\pm\sqrt{(\disc(f^n(x)-\alpha))}$ and fixes $\sqrt{\disc(f^m(x)-\alpha)}$ for all $0\leq m<n$. This proves the desired result.
\end{proof}
The above result indicates that we may overlook many quadratic extensions if the polynomial $f$ is PCF. As Proposition~\ref{prop:discriminant}, the discriminant of a polynomial is the product of the numbers on the orbits of its critical points. Hence, if a polynomial is PCF, then the square roots of the discriminants in Proposition~\ref{prop:discriminant} would be contained in some finite extension of the field of coefficients. This conclusion is further stated in the following lemma, which plays a crucial role in the proof of Theorem~\ref{intro-theorem-1}.

\begin{lemma}\label{lemma:periodicdisc}
Let $f$ be a PCF polynomial over $K$, and $\alpha\in K$. Let $L$ be the minimal integer such that $f^L(\mathcal{C}_f)$ is a periodic set, and let $O$ be the minimal integer such that $f^{L+O}(\mathcal{C}_f)=f^L(\mathcal{C}_f)$. Then, we have the following: 
\begin{enumerate}
    \item\label{lemma:periodicdisc_1} Suppose that all critical points are periodic, and the degree of $f$ is odd. Then $\disc(f^{2O}-\alpha)$ is a square in $K$. 
    \item\label{lemma:periodicdisc_2} Suppose that all critical points are periodic, and the degree of $f$ is even. Then, $\disc(f^{O+1}-\alpha)$ is a square in $K_f^1(\alpha)$.
    \item Suppose that a critical point is strictly preperiodic, and the degree of $f$ is odd. Then, $\disc(f^{L+2O-1}-\alpha)$ is a square in $K_f^{L-1}(\alpha)$. In particular, $K_f^{L-1}(\alpha)= K$ if $L=1$.
    \item Suppose that a critical point is strictly preperiodic, and the degree of $f$ is even. Then,  $\disc(f^{L+O}-\alpha)$ is a square in $K_f^{L}$.
\end{enumerate}
\end{lemma}
\begin{proof}
Let us define a "potential non-square factor" of a product $a=a_1^{e_1}\cdots a_n^{e_n}$ as $a_1^{e_1'}\cdots a_n^{e_n'}$, where $e_i'=0$ for even $e_i$ and $e_i'=1$ for odd $e_i$. This definition is ambiguous because it depends on the representation of $a$ as a product. However, it is still a useful concept. If the potential non-square factor of a product $a$ is a square, then $a$ is also a square.

We claim that the potential non-square factor of the product~\eqref{eq:disc} is
\begin{equation}\label{odd}
\left(\left(\dfrac{d}{4}\right)d\right)^{b_n}\left(\prod_{k=1}^n\prod_{c\in\mathcal{C}_f}(f^k(c)-\alpha)\right)
\end{equation}
for odd $d$, and is
\begin{equation}\label{even}
a_f\left(\prod_{c\in\mathcal{C}_f}(f^n(c)-\alpha)\right)
\end{equation}
for even $d$. We will prove this claim eventually, but let us discuss the consequence. We omit $c\in \mathcal{C}_f$ in the following product.
\begin{enumerate}
    \item We use Equation~\eqref{odd}. By plugging $2O$ for $n$, we see that the potential non-square factor is a square.
    \item We observe that $a_f\prod (f^{O+1}(c)-\alpha)=a_f\prod (f(c)-\alpha)$, and it implies that the potential non-square factor of $\disc(f^{O+1}-\alpha)$ is a square in $K_f^1(\alpha)$. 
    \item Note that each of the followings $f^L(\mathcal{C}_f),\ldots, f^{L+2O-1}(\mathcal{C}_f)$ repeat twice. Thus, the product
    \[
    \prod_{k=L}^{L+2O-1}\prod_{c\in\mathcal{C}_f}(f^k(c)-\alpha)
    \]
    is a square. Moreover both $L-1$ and $L+2O-1$ have the same parity, so we observe Equation~\eqref{odd} by plugging $L+2O-1$ for $n$ and conclude that $\disc(f^{L+2O-1}-\alpha)$ is a square in $K_f^{L-1}(\alpha)$.
    \item We have that $a_f\prod(f^L(c)-\alpha)=a_f\prod(f^{L+O}(c)-\alpha$), so $\disc(f^{L+O}-\alpha)$ is a square in $K_f^{L}(\alpha)$.
\end{enumerate}

Let us now clean up the remaining assertion. To do this, we will use induction. Consider the case when $d$ is odd. For the base case, we have
\[
\disc(f-\alpha)=\left(\dfrac{d}{4}\right)a_f^{d-1}d^d\prod_{c\in\mathcal{C}_f}(f(c)-\alpha),
\]
and it can be observed that this holds.

Assuming the statement is true for all integers less than $n$, the potential non-square part of
\begin{align}\label{eq: f^n-a}
    \disc(f^n-\alpha) & =\left(\dfrac{d}{4}\right)a_f^{d^{2n-1}-1}d^{d^n}(\disc(f^{n-1}-\alpha))^d\prod_{c\in\mathcal{C}_f}(f^n(c)-\alpha)
\end{align}
is given by
\begin{align*}
&\left(\dfrac{d}{4}\right)d\disc(f^{n-1}-\alpha)\prod_{c\in\mathcal{C}_f}(f^n(c)-\alpha)\\
=&\left(\dfrac{d}{4}\right)d\left[\left(\left(\dfrac{d}{4}\right)d\right)^{b_{n-1}}\prod_{k=1}^{n-1}\prod_{c\in \mathcal{C}_f}(f^k(c_i)-\alpha)\right]\prod_{c\in\mathcal{C}_f}(f^n(c)-\alpha).
\end{align*}
Therefore, the statement holds for odd $d$. 

For even $d$, we can check Equation~\ref{even} directly from Equation~\ref{eq: f^n-a}. Note that $\left(\dfrac{d}{4}\right)=1$ for even $d$. Furthermore, $d^d$ is always a square, and the power of $a_f$ is always odd. Thus, the potential non-square part of $\disc(f^n-\alpha)$ is
\[
a_f\prod_{c\in\mathcal{C}_f}(f^n(c)-\alpha).
\]
\end{proof}


\section{Main results}\label{sec:3}
\subsection{Universal embedding of dynamical groups of PCF maps}
The quadratic extensions are lost only after some finite extension. Equivalently, the PCF condition forces there
to be only finitely many quadratic extensions over the ground field. Following this idea, we can show the following well-known result by Jones (see~\cite{Rafe2013}).
\begin{corollary}
If $f$ is a PCF polynomial over a number field $K$, then there exists some integer $N$ such that $\Gal_f^N(\alpha)$ is not isomorphic to $[S_d]^N$. Furthermore, the index of the arboreal representation of $\Gal^\infty_f(\alpha)$ in $[S_d]^\infty$ is infinite
\end{corollary}
\begin{proof}
This directly follows from Theorem~\ref{intro-theorem-1} and Proposition~\ref{prop:unboundedindex}.
\end{proof}

Now, we are ready to proof Theorem~\ref{intro-theorem-1}.
\begin{proof}[Proof of Theorem~\ref{intro-theorem-1}]
Since $\alpha$ is not periodic, the dynamical tree of $f$ at $\alpha$ is isomorphic to a full $d$-ary tree, where $d$ is the degree of $f$. We will proceed by induction. Our induction hypothesis is that $\Gal_f^n(\alpha)$ is isomorphic to a subgroup of $E_n^{(m,m')}$ or $F_n^{(m, m')}$ for all $n\geq 1$. This hypothesis is trivial for $n<m$ because we set $E_n^{(m,m')}$ and $F_n^{(m,m')}$ to be $\Aut(T_n)$ for $n<m$. Thus, our base case starts from $n=m$.

First, assume $L\leq 1$ and the degree of $f$ is odd, and let $m=2O$. By Lemma~\ref{lemma:periodicdisc}.(\ref{lemma:periodicdisc_1}), $\disc(f^{m}-\alpha)$ is a square in $K$, so the dynamical Galois group $\Gal_f^{m}(\alpha)$ must be isomorphic to a subgroup of $\Aut(T_m)\cap \ker(\sgn_m)$. This proves the initial step $n=m$. Now, assuming $\Gal_f^m(\alpha)$ is isomorphic to a subgroup of $E_n^m$ for $n>m$, we shall prove the statement for $n+1$. Let $f^{-1}(\alpha)=\{y_1,\ldots, y_d\}$. By the induction hypothesis, the dynamical Galois group $\Gal(K_f^{n}(y_i)/K(y_i))$ is isomorphic to a subgroup of $E_n^m$. Clearly, $\Gal(K_f^{n+1}(\alpha)/K)$ is isomorphic to a subgroup of $E_n^m\wr \Aut(T_1)$. Since the discriminant of $f^{m}-\alpha$ is a square in $K$, any permutation acting on the $n+1$-th level of $T_n$ that is restricted to the $m$-th level must be even. Hence, $\Gal_f^{n+1}(\alpha)$ is contained in the kernel of $\sgn_m$. This shows that $\Gal_f^{n+1}(\alpha)$ is isomorphic to a subgroup of $E_m^{n+1}$, completing the induction step.

Now suppose $f$ is odd and $L>1$. Then, by Lemma~\ref{lemma:periodicdisc}.(\ref{lemma:periodicdisc_2}), when $m=L+2O+1$ and $m'=L-1$, we have that $\disc(f^{m}-\alpha)$ is a square in $K_f^{m'}(\alpha)$. This means that any element $\sigma\in \Gal_{\text{top}}\coloneqq\Gal(K_f^{m}(\alpha)/K_f^{m'}(\alpha))$ is an even permutation. For the initial step $n=m$, we consider the following diagram:
\[
\begin{tikzcd}
1\arrow[r]&\Gal_{\text{top}}\arrow[r]\arrow[d,"\sgn_m"]&\Gal_f^{m}(\alpha)\arrow[r,"res_{m'}"]\arrow[d,"\sgn_m"]&\Gal_f^{m'}(\alpha)\arrow[r]\arrow[d,"\sgn_{m'}"]&1\\
1\arrow[r]&1\arrow[r]&C_2\arrow[r,"id"]&C_2\arrow[r]&1
\end{tikzcd}
\]
where the homomorphism from $\Gal_{\text{top}}$ to $\Gal_f^m(\alpha)$ is the canonical inclusion. We observe that both $\sqrt{\disc(f^{m}-\alpha)}$ and $\sqrt{\disc(f^{m'}-\alpha)}$ have the same potential non-square factor over $K$. This means that for $\sigma\in \Gal_f^{m}(\alpha)$, both $\sigma$ and $\res_{m'}(\sigma)$ must have the same parity. Consequently, the bottom line of the diagram is exact and the diagram is commutative—that is, $\sgn_m(\sigma)=\sgn_{m'}(\sigma)$ for all $\sigma\in \Gal_f^{m}(\alpha)$. By Lemma~\ref{lemma:sgn2}, since $\sigma\in\ker(\sgn_i^{(m,m')})$, $\Gal_f^{m}(\alpha)$ is isomorphic to a subgroup of $E_{m}^{(m,m')}$. This proves the base case when $n=m$. For the induction step, an argument similar to the case when $L\leq 1$ and $\deg f$ is odd shows that $\Gal_{f}^{n+1}(\alpha)$ is isomorphic to a subgroup of $E_{n+1}^{(m,m')}$.

For the remaining cases, we simply choose $m$ and $m'$ according to Lemma~\ref{lemma:periodicdisc}, and all arguments remain the same.
\end{proof}
\begin{remark}
The finding of a square discriminant is one way to demonstrate an infinite Odoni index. If the Galois group $G_f^1(\alpha)$ is a nontrivial subgroup of $S_d$, then the Odoni index of $f$ is infinite. A less explored family is the following: consider all the critical points $c_1,\ldots, c_{2r}$ of $f$, including multiple roots as distinct points. If there exists an $n_i$ such that $f^{n_i}(c_{2i-1})=f^{n_i}(c_{2i})$ for each $i=1,\ldots r$, then the discriminant of $f^n-\alpha$ is a rational square for sufficiently large values of $n$. This family has not been thoroughly studied yet.
\end{remark}

\subsection{The rank of \texorpdfstring{$E_n^2$}{LG}}\label{sec:the rank of E_n}
Let $S$ be a subset of a group $G$, iff $G=\langle S\rangle$, we call $S$ a \textit{generating set} of $G$. It is a long-term goal to determine the cardinality of \textit{minimal generating set} of $G$. Let
\[
d(G)=\min\{\# S\mid G=\la S\ra\}.
\]
Here we recall an important theorem regarding $d(G)$ by Dalla Volta and Lucchini (see~\cite{DallaVoltaLucchini}).
\begin{lemma}\label{lemma:boundedrank}
If a finite non-cyclic group $G$ contains a unique minimal normal
subgroup $M$, then we have
\[
d(G)\leq\max\{2,d(G/M)\}.
\]
\end{lemma}
A consequence of this lemma is that if there exists a tower of normal subgroups as follows:
\[
\{1\}= N_1\triangleleft N_2\triangleleft\cdots\triangleleft N_k\triangleleft G \label{eq:towerofnormalsubgroup} \tag{\P}
\]
such that each factor $N_{i+1}/N_i$ is the unique minimal normal subgroup in the quotient $G/N_i$, then $d(G)\leq \max\{2,d(G/N_i)\}$ for all $i$.

In the following, we want to show that $E_n^2$ has a tower of normal subgroups satisfying Condition~(\eqref{eq:towerofnormalsubgroup}). We apply Argument~\ref{two-step} to show that a proper normal subgroup $N\triangleleft G$ is the unique minimal one.
\begin{argument}[H]
\centering
\begin{minipage}{.9\linewidth}
\begin{enumerate}
    \item For any proper normal subgroup $N'$ of $G$, we have $N\cap N'\neq\{1\}$.
    \item $N$ is a minimal normal subgroup.
\end{enumerate}
\end{minipage}
\caption{Two-step argument}
\label{two-step}
\end{argument}
\begin{lemma}\label{lemma:scsneqc}
\begin{enumerate}
    \item For a nontrivial element $\sigma\in E_n^m$, there exists an element $\mathbf{c}\in \ker(\res_{n-1}:E_{n}^m\to E_{n-1}^m)$ such that 
    \[
    \sigma \mathbf{c}\sigma^{-1}\mathbf{c}^{-1}\neq 1.
    \]
    \item\label{lemma:scsneqc.b} Let $\mathbf{c}\in \ker(\res_{n-1}:E_{n}^m\to E_{n-1}^m)$. If the order of $\mathbf{c}$ is $2$, then there exists an element $\sigma\in E_n^m$ such that 
    \[
    \sigma\mathbf{c}\sigma^{-1}\mathbf{c}\in M\setminus\{1\},
    \]
    where $M=\{((a_i)_{i\in I};1)\mid a_i\in A_d\}$.
\end{enumerate}
\end{lemma}
\begin{proof}
\begin{enumerate}
    \item For any nontrivial element $\sigma=((a_i)_{i\in I};b)\in E_n^m\subset \Aut(T_1)\wr\Aut(T_{n-1})$, we can write $\sigma$ as the product of $\mathbf{b}$ and $\mathbf{a}$ where $\mathbf{a}=((a_i)_{i\in I};1)$ and $\mathbf{b}=((1)_{i\in I}; b)$. Hence, we can write $\sigma\mathbf{c}\sigma^{-1}$ as
    \[
    \mathbf{b}\mathbf{a}\mathbf{c}\mathbf{a}^{-1}\mathbf{b}^{-1}.
    \]
    To prove the first statement, we aim at find an element $\mathbf{c}$ such that $\sigma\mathbf{c}\sigma^{-1}\neq \mathbf{c}$.
    
    If $\mathbf{a}$ is the identity, we only need to choose $\mathbf{c}={((c_i)_{i\in I};1)}$ such that $c_i\neq c_j$ for some indices $i$ and $j$ with $\mathbf{b}^{-1}(i)=j$.
    
    If both $\mathbf{a}$ and $\mathbf{b}$ are not identity with indices $i$ and $j$ satisfying $\mathbf{b}^{-1}(i)=j$, then, for $d\geq 5$, we can choose $c_i$ and $c_j$ that are products of different length of disjoint cycles respectively, e.g., $c_i$ and $c_j$ can be a product of two disjoint $2$-cycles and a 3-cycle respectively. For $d=3$, we may choose $c_1$ to be a 3-cycle and choose $c_2$ and $c_3$ to be a $2$-cycle. Then $a_ic_ia_i^{-1}\neq a_jc_ja_j^{-1}$ due to our  choice of permutation. Treating $\mathbf{a}\mathbf{c}\mathbf{a}^{-1}$ as a new $\mathbf{c}$, we revert to the previous case. 
    
    Finally, if $\mathbf{b}$ is the identity, the statement is trivially true because $\ker(\res_{n-1}:E^m_n\to E^m_{n-1})$ is non-commutative.
    \item Note that if the order of $\mathbf{c}$ is $2$, then $\mathbf{c}^{-1}=\mathbf{c}$. Write $\mathbf{c}=((c_i)_{i\in I};1)$, we can definitely find some $a_1\in S_d$ such that $a_1c_1a_1^{-1}\neq c_1$. Thus, we will arbitrary choose $a_i$ for $i\neq 1$ such that $\mathbf{a}=((a_i)_{i\in I};1)\in E_n^m$, and so $\mathbf{a}\mathbf{c}\mathbf{a}^{-1}\mathbf{c}$ is not identity. In the end, we clearly has even permutation on each coordinate of $\mathbf{a}\mathbf{c}\mathbf{a}^{-1}\mathbf{c}$, so it belongs to $M$.
\end{enumerate}
\end{proof}
\begin{lemma}\label{lemma:inverseautomorphismisbyconjugate}
Let $\sigma\in [S_d]^n$. There exists an $\tau\in [S_d]^n$ such that $\sigma^\tau=\sigma^{-1}$.
\end{lemma}
\begin{proof}
We proceed by induction. For the base case, the statement trivially follows from the fact that $\Aut(S_d)\cong S_d$. 

We assume the statement is true for a positive integer $N$. Given an element $\sigma\in \Aut(T_{N+1})$, we may write $\sigma$ as
\[
\sigma=((\sigma_i)_{i\in I(T_N)};\sigma_0)\in\Aut(T_1)\wr\Aut(T_{N})
\]
for $\sigma_i\in S_d$ and $\sigma_0\in\Aut(T_N)$. Thus, for each $i$, we have $\tau_i$ such that $\sigma_i^{\tau_i}=\sigma_i^{-1}$ by the fact that $\Aut(S_d)\cong S_d$ and the induction hypothesis. Now, we observe
\[
((\tau_{\tau_0^{-1}(i)})_{i\in I(T_d)};1)(1;\tau_0)((\sigma_i)_{i\in I(T_d)};\sigma_0)(1;\tau_0)^{-1}((\tau_{\tau_0^{-1}(i)})_{i\in I(T_d)};1)^{-1}
\]
is the inverse of $\sigma$.
\end{proof}
\begin{lemma}\label{lemma:elementinres}
Suppose one of the following holds:
\begin{enumerate}
    \item Let $G$ be a group. Assume that, for any $g\in G$, there exists $g'\in G$ such that $g^{g'}\neq g^{-1}$.
    \item Let $\res_n:G\wr\Aut(T_n)\to \Aut(T_n)$ be the restriction map. Assume $H\subseteq \ker(\res_n)$ is a normal subgroup of $G\wr\Aut(T_n)$, and that, for any $h\in H$, there exists $\sigma\in\Aut(T_n)$ such that $h^\sigma \neq h^{-1}$.
\end{enumerate}
 Then, any proper normal subgroup of $G\wr\Aut(T_n)$ or $H\ltimes\Aut(T_n)$ contains a nontrivial element in $\{((g_1,\ldots,g_n);1)\mid g_i\in G\}$
\end{lemma}
\begin{proof}
For any element $\tau\in N\triangleleft G\wr\Aut(T_n)$, there exist $\mathbf{g}=((g_i)_{i\in I(T_n)};1)$ and $\mathbf{s}=(1;\sigma)$ for $g_i\in G$ and $\sigma\in \Aut(T_n)$ such that $\tau=\mathbf{s}\mathbf{g}$. By Lemma~\ref{lemma:inverseautomorphismisbyconjugate}, there exists $\mathbf{t}\in G\wr\Aut(T_n)$ such that $\mathbf{t}\mathbf{s}\mathbf{t}^{-1}=\mathbf{s}^{-1}$. Thus, we have
\[
\tau^{\mathbf{t}}\tau=((\ast);1),
\]
where $\ast$ means an arbitrary element.

To prove that the first condition implies the conclusion, one can conjugate some $g_i\in G$ by $g_i'\in G$ such that $g_i^{g_i'}\neq g_i^{-1}$. Thus, one may conjugate $\tau$ by some element $\mathbf{g}'\in \ker(\res_n)$ such that $\tau^{\mathbf{t}}\tau^{\mathbf{g}'}\neq 1$.

Next, we prove the second condition implies the conclusion. Since $H$ is normal in $G^{d^n}$, the semidirect product $H\ltimes\Aut(T_n)$ is well-defined. Following a similar manner as before, we can construct a nontrivial element in $\ker(\res_n)$.
\end{proof}
Let $e_i$ be an element of the form $(1,\ldots,-1,\ldots,1)$ in $C_2^{d^{n-1}}$, where $-1$ appears only in the $i$-th coordinate. Let $k_i=d^{i-1}+1$. We define
\[
X_i=\{(e_1+e_{k_i};1)^{\sigma}\mid \sigma\in \Aut(T_n)\}\subseteq C_2^{d^n}\ltimes \Aut(T_{n})
\]
for $i=1,2,\ldots, n$. The following lemma is crucial.
\begin{lemma}\label{lemma:C2AutTn}
Let $H_i=\la X_i\ra\ltimes\Aut(T_n)$. Then, $H_i=H_1$ for all $i$, and $\la X_1\ra$ is the unique minimal normal subgroup of $H_i$. In particular, for any $i$ and $n$, let $\res_n:H_i\to\Aut(T_n)$, then we have $\langle X_1\rangle=\ker(\res_n)$.
\end{lemma}
\begin{proof}
To show $\langle X_1\rangle$ is the unique minimal proper subgroups of $H_i$, we apply Argument~\ref{two-step}. Firstly, we notice that every element in $X_1$ has order $2$. Secondly, $H_n$ is a subgroup of $C_2\wr\Aut(T_n)$, and, for any $h\in \la X_i\ra$, we can find $\sigma\in \Aut(T_n)$ such that $h^\sigma\neq h$. By Lemma~\ref{lemma:elementinres}, the intersection of any proper normal subgroup of $H_n$ and the kernel of $\res:H_n\to\Aut(T_n)$ contains a nontrivial element $h$. We will use this nontrivial element to construct $e_1+e_2$ by a sequence of conjugations. 

Since $\Aut(T_{n})$ acts transitively on $I(T_n)$, we can assume that the nontrivial element $h$ in the intersection of a proper normal subgroup of $H_n$ and the kernel of $\res_n:H_n\to\Aut(T_n)$ has $-1$, the nontrivial element of $C_2$, in the block of $\mathcal{B}(T_n/T_{n-1})(1,\ldots, 1)$. Moreover, $\Aut(T_n)$ acts by conjugation, so $\sigma=(((12\cdots d),1,\ldots,1);1)$ shifts the block. Thus, $h^\sigma h$ has even many $-1$ in the first block and $1$ on all other entries. Next, we claim that, for even many $-1$ and odd many $1$ in a black with $S_d$ acting on the indices, we can always obtain
\[
(-1,\ldots, -1,1,\ldots, 1)
\]
by conjugating an even permutation. If $h$ is already in the desired form, then we are done. Otherwise, there are three consecutive indices, say $i$, $i+1$, and $i+2$, with $1$ and $-1$ distributed in one of the following layouts
\begin{align*}
    1, -1, -1;\\
    -1, 1, -1;\\
    1, -1, 1;\\
    1, 1, -1.
\end{align*}
Applying a 3-cycle $(i, i+1, i+2)$ several times by conjugate, one can move the above layouts to either $-1,-1,1$ or $-1,1,1$. Since there are only a finite number of indices, the desired form can be achieved by applying a finite number of conjugations. Since we apply $3-cycle$ each time, the desired result can be achieved by an even permutation. Now, we can conjugate $(-1,-1,\ldots, -1,1,\ldots,1)$ by $(12\cdots d)$ to get $(1,-1,\ldots, -1,1,\ldots, 1)$. Multiplying $(-1,\ldots,-1,1,\ldots,1)$ and $(1,-1,,\ldots,-1,1,\ldots,1)$ gives $(-1,1,\ldots,1,-1,1,\ldots,1)$. Finally, we repeating the process of conjugating a 3-cycle $(i,i+1,i+2)$ to get $(-1,-1,1,\ldots, 1)$. This shows that an arbitrary proper normal subgroup contains the generators of $H_1$, which implies that $\la X_1\ra$ is the unique minimal normal subgroup of $H_i$.
\end{proof}
\begin{remark}
The above proof can be written as the following psudocode.

\begin{algorithm}[H]
\KwData{A tuple $\mathbf{a}=(a_1,\ldots, a_d)\in C_2^d$ with $\prod{a_i}=1$ and odd $d$}
\KwResult{The tuple $\mathbf{a}=(-1,-1,1,\ldots,1)$}
\BlankLine
$\mathbf{a}'\coloneqq\mathbf{a}$\;
\While{$\mathbf{a}\neq(-1,-1,1,\ldots1)$}{
    $\mathbf{a}\coloneqq\mathbf{a'}$\;
    $\mathbf{b}\coloneqq(-1,\cdots,-1,1,\cdots,1)$  with same number of $-1$ as in $\mathbf{a}$\;
    \While{$\mathbf{a}\neq\mathbf{b}$}{
        Find the minimal index $i$ such that $a_i==1$ and $a_{i+1}==-1$\;
        $(a_i,a_{i+1},a_{i+2})\coloneqq (a_{i+1},a_{i+2},a_{i})$\;
    }
    $\mathbf{a}'\coloneqq(a_d',a_1',a_2',\ldots,a_{d-1}')$\;
    $\mathbf{a}'\coloneqq\mathbf{a}\mathbf{a}'$\;
}
\caption{The pseudocode of the above proof}
\end{algorithm}
\end{remark}
\begin{lemma}\label{lemma:a1a2a3a4}
For $\mathbf{a}=((a_1,a_2,a_3,\ldots, a_d);1)\in E_2^2$, there exists $\sigma\in E_2^2$ such that $\mathbf{a}^{\sigma}=((a_1,a_2^{-1},a_3^{-1},\cdots,a_d^{-1});1)$ or $((a_1,a_2,a_3^{-1},\cdots,a_d^{-1});1)$.
\end{lemma}
\begin{proof}
For any $a_i$, there exists an element $s_i\in S_d$ such that $a_i^{s_i}=a_i^{-1}$. If $((1,s_2,\ldots, s_d);1)$ is in $E_2^2$, then the proof is complete. Otherwise, there are an odd number of odd permutations $s_i$. Without lose of generality, we can assume $s_2$ is odd. Then, $((1,1,s_3,\ldots, s_{d});1)$ is in $E_2^2$, and we can conjugate $\mathbf{a}$ by this element to get $((a_1,a_2,a_3^{-1},\cdots,a_d^{-1});1)$. 
\end{proof}

\begin{lemma}\label{lemma:E_n->E_n-1}
Let $d$ be an odd integer, and let $\res_{n}:E^2_{n+1}\to E^2_{n}$ be the natural restriction from $\Aut(T_{n+1})$ to $\Aut(T_n)$. Let $M_n=((A_d)_{i\in I},1)\subseteq E_{n+1}\subseteq \Aut(T_1)\wr\Aut(T_n)$ where $A_d$ is the alternative group of degree $d$. Then, $E_{n}^2$ has a tower of normal subgroups
\[
\{1\}\triangleleft M_n\triangleleft\ker(\res_{n})\triangleleft E_{n+1}^2
\]
satisfying Condition~(\ref{eq:towerofnormalsubgroup}).
\end{lemma}
\begin{proof}
We first show that $M_n$ is a normal subgroup. It is clear that $M_n$ is a proper subgroup of $E_{n}^2$, and, for an element $\sigma\in E_{n}^2$, we can express $\sigma$ as $\sigma=\mathbf{b}\mathbf{a}$ where $\mathbf{a}=((a_i)_{i\in I},1)$ and $\mathbf{b}=((1)_{i\in I},b)$. Therefore, we have
\[
\sigma M_n\sigma^{-1}=\mathbf{b}\mathbf{a}M_n\mathbf{a}^{-1}\mathbf{b}^{-1}=\mathbf{b}M_n\mathbf{b}^{-1}=M_n,
\]
which implies that $M_n$ is a normal subgroup of $E_n^2$.

Next, we show that a nontrivial normal subgroup $N$ contains an element in $M_n$. Given an element $\sigma\in N$ and $\mathbf{c}=((c_i)_{i\in I},1)\in E^2_{n}$, we have $\sigma\mathbf{c}\sigma^{-1}\mathbf{c}^{-1}\in N$. If $\sigma$ is in $M_n$, then we are done. If not, then we use Lemma~\ref{lemma:scsneqc} to show that for a given $\sigma$, there exists $\mathbf{c}\in \ker(\res_{n-1})$ such that
\[
\sigma\mathbf{c}\sigma^{-1}\mathbf{c}^{-1}\neq 1.
\]
Since $\mathbf{c}$ is also in a normal subgroup, it follows that $1\neq \sigma\mathbf{c}\sigma^{-1}\mathbf{c}^{-1}\in N\cap \ker(\res_n)$. Let $\gamma = \sigma\mathbf{c}\sigma^{-1}\mathbf{c}^{-1}$. If $\gamma$ has order greater than $2$, then $\gamma^2\in N\cap M_n\setminus\{1\}$. If not, by Lemma~\ref{lemma:scsneqc}.\ref{lemma:scsneqc.b}, we find $\delta\in E_n^2$ such that $\delta\gamma\delta^{-1}\gamma^{-1}\in N\cap M_n\setminus\{1\}$.

If we can prove that $M_n$ is a minimal normal subgroup of $E^2_{n}$, then we can complete Argument~\ref{two-step} and conclude that $M_n$ is the unique normal subgroup of $E^2_{n}$. Using induction, we show that $M_n$ is a minimal normal subgroup of $E^2_n$ for $n\geq 2$.

For the base case $E^2_2$, we define a homomorphism $\rho_i$ from $\ker(\res_1:E_2^2\to E_1^2)$ to $\Aut(T_1)$ as follows:
\begin{align*}
    \rho_i((a_i)_{i\in I};1)= a_i
\end{align*}
If $N$ is a nontrivial normal subgroup of $E_2^2$, then $M_2\cap N$ is a normal subgroup of $\ker(\res_1)$. Hence, $\rho_i(M_2\cap N)$ is a normal subgroup of $A_d=\rho_i(m')$. Since $A_d$ and $\{1\}$ are the only normal subgroups of $A_d$ and $\res_{1}(E_2^2)=\Aut(T_1)$ acts transitively on the level $1$ of the tree, we conclude that $\rho_i(M_2\cap N)\cong A_d$ for all $i$. Without loss of generality, let $i=1$. Then there exists a nontrivial element $\mathrm{a}=((a_1,a_2,\ldots, a_{d});1)\in N\cap M_2$ for some $a_1\in A_d$. We may choose $a_1=(123)$ or $(12\cdots n)$, as they form a generating set of $A_d$ for odd $d$.  We claim that $\mathrm{a}'=((a_1,1,\ldots, 1);1)$ exists in the intersection. By Lemma~\ref{lemma:a1a2a3a4}, we know that we can conjugate $\mathrm{a}$ to either $\mathrm{a}_1=((a_1,a_2^{-1},\ldots, a_d^{-1});1)$ or $\mathrm{a}_2=((a_1,a_2,a_3^{-1},\ldots, a_d^{-1});1)$. If the conjugation results in $\mathrm{a}_1$, then a power of $\mathrm{a}\mathrm{a}_1$ is in the desired form. If it results in $\mathrm{a}_2$, then we can conjugate $\mathrm{a}\mathrm{a}_2$ by either $((1,s,1,1,\ldots, 1);1)$ or $((1,s,(1,2),1,\ldots,1);1)$ for some $s\in S_d$, where one of these permutations is in $E_2^2$. This yields $s\mathrm{a}\mathrm{a}_2s^{-1}=((a_1^2,a_2^{-2},1,\ldots,1);1)$. Hence, the desired form can be achieved by taking a power of $\mathrm{a}\mathrm{a}_2s\mathrm{a}\mathrm{a}_2s^{-1}$. We therefore prove that $M_2$ is a minimal normal subgroup of $E_2^2$. 

For the inductive step, let us assume that $M_n$ is a minimal normal subgroup of $E_n^2$. This implies 
\[
\mathcal{A}=\{((a,1,\ldots,1);1)\in\Aut(T_1)\wr\Aut(T_{n-1})\mid a\in A_d\}
\]
is a subgroup of $M_n\subseteq E_n^2$, and thus
\begin{align*}
\mathcal{A}'&=\{((a',1,\ldots,1);1)\in E_n^2\wr\Aut(T_1))\mid a'\in \mathcal{A}\}\\
&=\{((a,1,\ldots,1);1)\in \Aut(T_1)\wr\Aut(T_n))\mid a\in A_d\}
\end{align*}
is a subgroup of $M_{n+1}$. By Lemma~\ref{lemma:transitive}, it follows that $M_{n+1}$ is a minimal normal subgroup of $E_{n+1}^2$. Hence, we finish the induction step and can conclude that $M_{n}$ is the minimal normal subgroup of $E_n^2$ for $n\geq 2$.

We claim that $E_n^2/M_n$ is isomorphic to $C_2\wr\Aut(T_{n-1})$. To demonstrate this, we define a homomorphism
\begin{align*}
    \phi:&E_n^2\subseteq\Aut(T_1)\wr\Aut(T_{n-1})\to C_2\wr\Aut(T_{n-1})\\
    &((a_i)_{i\in I(T_{n-1})},b)\mapsto ((\sgn(a_i))_{i\in I(T_{n-1})},b).
\end{align*}
It is straightforward to verify that $\ker(\phi)\cong M_n$ and that $\phi$ is surjective. Finally, by applying Lemma~\ref{lemma:C2AutTn}, we conclude that $\ker(\res_{n-1})/M_n$ is the minimal normal subgroup of $E_n^2/M_n$, which is the final assertion of this lemma.
\end{proof}
We are now ready to prove Theorem~\ref{intro-theorem-2}.
\begin{proof}[Proof of Theorem~\ref{intro-theorem-2}]
By Lemma~\ref{lemma:boundedrank} and Lemma~\ref{lemma:E_n->E_n-1}, we have 
\[
d(E_n^2)\leq \max\{2,d(E_n^2/M_n)\}\text{ and }d(E_n^2/M_n)\leq \max\{2,d(E_n^2/\ker(\res_n))\},
\]
which implies $d(E_n^2)\leq \max\{2,E_n^2/\ker(\res_n)\}$. As $E_n^2/\ker(\res_n)\cong E_{n-1}^2$, we have
\[
d(E_n^2)\leq\max\{2,d(E_{n-1}^2)\}\leq \cdots\leq\max\{2, d(E_1^2)\}=2.
\]

To prove $E_n^2$ has a unique chief series, we argue by induction. When $n=1$, we see that $S_d$ has a unique chief series. Since $d$ is odd, we have
\[
1\triangleleft A_d\triangleleft S_d.
\]
Let us assume $E_n$ has a unique chief series for some $n$. Now, given a nontrivial normal subgroup $N$ in $E_{n+1}$, since $M_{n+1}$ is the unique minimal normal subgroup of $E_{n+1}$, we must have $N\supseteq M_{n+1}$. Therefore, $N/M_{n+1}$ is either a trivial normal subgroup of $E_{n+1}/M_{n+1}$ or a normal subgroup contains $\ker(res_n)/M_{n+1}$. This shows that any nontrivial normal subgroup $N$ of $E_{n+1}$ contains $\ker(\res_n)$ or is equal to $M_{n+1}$. This shows that a normal series of $E_n$ lifts uniquely to a normal series of $E_{n+1}$ and contains $\ker(\res_n)$ by taking inverse of $\res_n$. Because the chief series in $E_n$ is unique by induction hypothesis, the chief series in $E_{n+1}$ is also unique.
\end{proof}
\subsection{Examples}
In this section, we will examine the arboreal representations of the polynomial $f(z)=2z^3-3z^{2}+1$ at a carefully selected base point $x$. The method employed here closely mirrors the approach in \cite{B-F-H-J-Y-2017}.

Let $K$ be a number field and let $x\in K\setminus\{0,1\}$. We assume the pair $(K,x)$ satisfies the following condition:
\begin{equation}\label{Condition}
\begin{array}{c}
     \text{there exist primes $\mathrm{p}$ and $\mathrm{q}$ of $K$ lying above $2$ and $3$ such that}  \\
     \text{either $v_\mathrm{q}(x)=1$ or $v_\mathrm{q}(1-x)=1$, and either $v_\mathrm{p}(1-x)= 1$ or $v_\mathrm{p}(x)=1$}
\end{array}
\end{equation}
\begin{lemma}
If $(K,x)$ satisfies Condition~(\ref{Condition}), then $f^n(z)-x$ is Eisenstein at $\mathrm{q}$ for all $n\geq 1$. In particular $f^n(z)-x$ is irreducible.
\end{lemma}
\begin{proof}
For $n=1$ and $v_\mathrm{q}(x)=1$, we have
\[
f(z)\equiv (-z^3+1)\equiv -(z-1)^3\mod \mathrm{q},
\]
so $f(z-1)-x$ is Eisenstein. For $n=2$, we have
\[
f^2(z)=-(f(z)-1)^3\equiv -(-z^3)^3=z^{9}\mod \mathrm{q}.
\]
Then, it is easy to prove the statement by induction.

For $v_\mathrm{q}(1-x)=1$, we have
\[
f(z)\equiv(-z^3)\mod\mathrm{q}
\]
and
\[
f^2(z)\equiv -f(z)^3\equiv (z-1)^3\mod \mathrm{q}.
\]
Again, the conclusion will follow easily by induction.
\end{proof}
\begin{lemma}\label{lemma:ramification index}
Let $(K,x)$, $\mathrm{p}$ and $\mathrm{q}$ satisfy Condition~(\ref{Condition}), and let $y\in f^{-n}(x)$.  Then:
\begin{itemize}
    \item There are primes $\mathrm{p}'$ and $\mathrm{q}'$ of $K(y)$ lying above $\mathrm{p}$ and $\mathrm{q}$ respectively such that
    \[
    e(\mathrm{p}'/\mathrm{p})=2^n\text{ and }e(\mathrm{q}'/\mathrm{q})=3^n
    \]
    \item $(K(y),y)$ satisfies Condition~(\ref{Condition}).
\end{itemize}
\end{lemma}
\begin{proof}
    We first prove the assertion for $\mathrm{q'}$. Suppose $v_\mathrm{q}(x)=1$, then $y\in f^{-1}(x)$ has $(y-1)^3\mod q$ which implies there exists a unique prime $\mathbf{q}'$ of $K(y)$ above $q$ satisfies $v_{\mathrm{q}'}(1-y)=1$, i.e. $e(\mathrm{q}'/\mathrm{q})=3$. On the other hand, if we have $v_{\mathrm{q}}(1-x)=1$, then it implies there exists a unique prime $\mathrm{q}'$ of $K(y)$ such that $v_{\mathbf{q}'}(y)=1$. 
    
    For $\mathrm{p}'$, we also prove the statement case by case. Let $v_\mathrm{p}(1-x)=1$. We note that the Newton polygon of $f(z)-x$ has a segment of length $2$ and height $1$ which implies there exists a prime $\mathrm{p}'$ above $\mathrm{p}$ such that $e(\mathrm{p}'/\mathrm{p})=2$. Moreover we have $v_\mathrm{p'}(y)=1$ for $y\in f^{-1}(x)$. Otherwise, let $v_\mathrm{p}(x)=1$. We then have
    \[
    f(z)-x\equiv -z^2+1\equiv -(z-1)^2\mod\mathrm{p},
    \]
    so there exists a prime $\mathrm{p}'$ over $\mathrm{p}$ with $v_{\mathrm{p}'}(1-y)=1$.
    
    Since the above holds for any arbitrary choice of $x$ with corresponding $y\in f^{-1}(x)$ satisfying Condition~(\ref{Condition}), the statement follows by induction."
\end{proof}

\begin{corollary}
Let $(K,x)$ be as in Proposition~\ref{lemma:ramification index}. Then, the degree of the extension $K_f^n(x)$ over $K$ is a multiple of $6^n$. Specifically, there exist elements of order $2$ and $3$ in $\Gal(K_f^n/K_f^{n-1})$.
\end{corollary}
\begin{proof}
The proof follows by induction, and more details can be found in~\cite{B-F-H-J-Y-2017}. The elements of order $2$ and $3$ exist by Cauchy's Theorem.
\end{proof}

\begin{proof}[Proof of Theorem~\ref{intro-theorem-3}]
In~\cite{B-F-H-J-Y-2017}, the authors showed that the existence of elements of order $2$ and $3$ implies the presence of a two-cycle permutation in $\Gal(K_f^2/K)$ and a three-cycle permutation in $\Gal(K_f^3/K)$ that fixes $f^{-1}(x)$. The same proof works in our case.
\end{proof}


\section{Questions}\label{sec:4}
We close by posing some questions. Firstly, we consider the choice of the polynomial $f(z)=2z^3-3z^2+1$. Pink (\cite{pink2013profinite}) has analyzed all degree $2$ PCF polynomials, so we focus on degree $3$ polynomials. For degree $3$ PCF polynomials, we can always arrange the critical points as ${0,1,\infty}$ by conjugation. According to Lemma~\ref{lemma:periodicdisc}, the simplest case is when $f(\mathcal{C}_f)=\mathcal{C}_f$. A case that has not been considered is $f(0)=1$ and $f(1)=0$, leading to the polynomial defined above. A more general polynomial with the same critical points portrait can be defined as follows. Let $B_d(x)$ be a degree $d$ Belyi polynomial with critical points ${0,1,\infty}$. Then $-B_d(z)+1$ transposes $0$ and $1$ while fixing $\infty$. Theorem~\ref{intro-theorem-2} implies that the $n$-th dynamical subgroup of $-B_d(z)+1$ is a subgroup of either $E_n^2(d)$ or $E_n^{(2,1)}(d)$.

The question we seek to answer is: "Is there an isomorphism between $\Gal_{-B_d+1}^n(\alpha)$ and $E_n^2$ or $E_n^{(2,1)}$?"

\begin{question}\label{question:belyi}
Given a number field $K$ and a degree $d$ polynomial $f\in K[z]$ with critical points portrait
\[
c_1\to c_2,\ c_2\to c_1,\ c_3\to c_3. 
\]
do we have $\Gal_f^n(\alpha)\cong E_n^2$ for a base point $\alpha\in K$?
\end{question}
A possible approach to answering Question~\ref{question:belyi} is through the method presented in~\cite{BEK2020}. This involves three steps: showing that the Geometric Galois group $\Gal(B^n(z)-t/\overline{K}(t))$, where $t$ is transcendental over $K$, is isomorphic to $E_n^2$ or $E_n^{(2,1)}$; deducing that the arithmetic Galois group is also isomorphic; and finally, studying the ramification condition of $K_f^n(\alpha)$ over $K$ to determine the desired isomorphism with a specified base point $\alpha$.

An interesting observation made by R. Lodge raises a further question: the polynomial used in this approach is a twisted Belyi polynomial. Specifically, if $l=(1-z)\in \Aut(B_n)$, then $f=-B(z)+1$ is simply $l\circ B$, and $f\circ f=B\circ B$. This raises the question: "do we always have $\Gal_f^n\cong \Gal_g^n$ if $f^k=g^k$ or, more generally, if $h\circ f^k=g^k\circ h$ for some degree $>1$ polynomial $h$?" This latter functional equation is known as the semiconjugate relation. In~\cite{Mythesis}, we studied the arboreal representation and semiconjugate relation and attempted to define dynamical isogeny and formulate a dynamical isogeny theorem.

While the answer to Question~\ref{question:belyi} may have already been provided by Lodge's observation, it is not the main focus of this paper and has been left for future research.

The following question naturally arises:
\begin{question}
Is there a pair $(f,\alpha)$ where $f$ is a PCF polynomial $f$ over $K$ with $\alpha\in K$ such that $\Gal_f^n(\alpha)$ is isomorphic to a finite index subgroup of $F_n^{(m,m')}$?
\end{question}
To address this question, a proper family of PCF maps of odd degree with tail length greater than but not equal to 1 must be found. However, there is no obvious candidate, making this an open question for future research. Moreover, finding a PCF polynomial that is isomorphic to a non-trivial subgroup of $E_n^{(m,m')}$ would also be of interest. For example, the rabbit polynomial has image lying in $E_n^6$, but what is the index of this realization?

Lastly, it would be interesting to determine the chief series of $E_n^m$ for $m\geq 3$ and the chief series of $E_n^{m,m'}$ and $F_n^{(m,m')}$. However, generalizing the argument to these cases is challenging due to the following observation. As noted in Lemma~\ref{lemma:C2AutTn}, let $H_n=\langle X_n\rangle$ with $n\geq 1$. Then, the quotient group $H_n/\langle X_1\rangle$ corresponds to the following exact sequence:
\[
\{1\}\to\la X_1\ra\to \la X_n\ra\to H_n\to \Aut(T_n)\to \{1\}.
\]
It can be seen that $\ker(\res_{n-1}:\Aut(T_n)\to\Aut(T_{n-1}))$ acts trivially on $\langle X_n\rangle/\langle X_1\rangle$. Thus, we have:
\[
H_n/\la X_1\ra\cong (\la X_n\ra/\la X_1\ra\times \ker(\res_{n-1}))\ltimes \Aut(T_{n-1}). 
\]
Since both $\langle X_n\rangle/\langle X_1\rangle$ and $\ker(\res_{n-1})$ are normal subgroups and have trivial intersection, a unique tower of normal subgroups satisfying Condition~(\ref{eq:towerofnormalsubgroup}) cannot be formed. The situation for $F_n^{(m,m')}$ appears to be more complex, and further investigation is needed before a conclusion can be drawn. Here, we formalize our final question:
\begin{question}
What is the number of chief series in $E_n^m$ for $m\geq 3$? Can we bound the rank of $E_n^m$ from a purely group theoretic perspective? Is there a unique chief series in $F_n^{(m,m')}$ for small $m$ and $m'$?
\end{question}
\section*{Acknowledgments}
The author is grateful to Prof. Alexander Hulpke, Prof. Russell Lodge, and Prof. Julie Wang for their constructive discussions and feedback, as well as to Prof. Andrea Lucchini for his valuable comments on the rank of the groups. The author also expresses his gratitude to his Ph.D. advisor, Prof. Thomas Tucker, and his postdoc mentor, Prof. L.-C. Hsia, for their support and valuable suggestions. The author would also like to acknowledge Alexander Galarraga for his constructive criticism of the manuscript. Lastly, the author extends his heartfelt thanks to the anonymous referee for their insightful comments, which greatly contributed to the quality of this paper and its eventual publication.
\bibliographystyle{plain}
\bibliography{references}
\end{document}